\documentclass{amsart}
\usepackage[T1]{fontenc}
\usepackage{amsmath,amssymb,amsthm}
\usepackage{tikz-cd}
\usepackage[hidelinks]{hyperref}
\hypersetup{
  pdftitle={A Determinant Proof of the Sun--Yang--Zuo Conjecture},
  pdfauthor={Jinbang Yang}
}

\newtheorem{theorem}{Theorem}
\newtheorem{corollary}[theorem]{Corollary}
\newtheorem{lemma}[theorem]{Lemma}
\newtheorem*{remark}{Remark}
\numberwithin{equation}{section}
\newcommand{\PP}{\mathbb P}
\newcommand{\OO}{\mathcal O}

\begin{document}

\title{A Determinant Proof of the Sun--Yang--Zuo Conjecture}
\author{Jinbang Yang}
\date{}

\begin{abstract}
We give an elementary proof of the determinant identity in
\cite[Conjecture~A.1]{SYZ} over any field of odd characteristic $p$.
The proof computes the determinant of a single linear map by
dividing by two different monic polynomials.
By \cite[Proposition~A.9]{SYZ}, the identity gives another proof
of the Sun--Yang--Zuo conjecture.
The conjecture identifies the Higgs--de Rham self-map of
the four-pointed projective line with the map induced by
multiplication by $p$ on its elliptic double cover.
\end{abstract}

\maketitle
\section{Introduction}

Let $k$ be a field of characteristic $p>2$, and let
$\lambda\in k\setminus\{0,1\}$. Put $p=2m+1$, let $t$
be an indeterminate, and set $g(x)=(x-1)(x-\lambda)$.
Define the coefficients $\gamma_n$ by
\[
  x^m g(x)^m=\sum_{n=m}^{3m}\gamma_n x^n.
\]
For $1\le n\le p-2$, set
\[
  \delta_n(\lambda,t)
  =\frac{\lambda^p(1-t)-(\lambda^p-t)\lambda^n}{n}.
\]
Let $A=A(\lambda,t)$ and $B=B(\lambda,t)$ be the matrices
with entries
\begin{align*}
  A_{rj}
    &=\delta_{m+j-r}(\lambda,t)
    &&(0\le r,j<m),\\
  B_{rs}
    &=
    \begin{cases}
      \gamma_{m+r-s}, & r\ge s,\\
      t\gamma_{m+p+r-s}, & r<s
    \end{cases}
    &&(0\le r,s\le m).
\end{align*}

The following theorem gives a polynomial form of the
determinant identity conjectured by Sun, Yang and Zuo \cite[Conjecture~A.1]{SYZ}.

\begin{theorem}[Determinant identity]\label{thm_det}
For $k$, $\lambda$, $A$, and $B$ as above, set
\[
  C=\left(\frac{1}{m+j-r}\right)_{0\le r,j<m},
  \qquad c_p=(-1)^m\det C.
\]
Then $c_p\in\mathbb F_p^\times$ and
\begin{equation}\label{eq_det}
  \det A(\lambda,t)
    =c_p\,\lambda^{m^2}(\lambda-1)^{m^2}\det B(\lambda,t)
  \qquad\text{in }k[t].
\end{equation}
\end{theorem}

The two matrices come from formulas for two self-maps of the projective line. Suppose now that $k$ is
perfect, and set $X=\PP^1_k$ with marked divisor
$D=0+1+\lambda+\infty$. The moduli space of stable rank-two
graded logarithmic Higgs bundles of degree one on $(X,D)$
is identified with $X$ by the unique zero of the nonzero
Higgs component
$\OO(1)\to\OO\otimes\Omega_X^1(\log D)$
\cite[\S4.1 and~\S4.3]{SYZ}. Under this identification, the
twisted Higgs--de Rham flow induces a morphism
$\varphi_{\lambda,p}\colon X\to X$ \cite[\S4.2 and Appendix A.1]{SYZ}.
If $a$ denotes the affine coordinate of the Higgs zero,
the denominator in that formula involves the determinant
of $A_p(\lambda,a)=A(\lambda,a^p)$, see \cite[Apendix~A and Theorem~A.5]{SYZ}.

For the second map, let $E_\lambda$ be the elliptic curve
$y^2=x(x-1)(x-\lambda)$ with origin at infinity. The
$x$-coordinate gives the double cover
$\pi\colon E_\lambda\to X$ branched over $D$.
Multiplication by $p$ commutes with $P\mapsto-P$ and
therefore descends through $\pi$ to a self-map of $X$,
the \emph{Latt\`es map} associated with $[p]$. The
denominator in the formula for this map involves the
determinant of $B_{m+1}(\lambda,a)=B(\lambda,a^p)$
\cite[Theorem~A.8]{SYZ}.

Sun, Yang and Zuo conjectured that these two self-maps
coincide \cite[Conjecture~4.8]{SYZ}.

Specializing \eqref{eq_det} at $t=a^p$ proves the determinant
identity of \cite[Conjecture~A.1]{SYZ}. 
\cite[Proposition~A.9]{SYZ} then yields \cite[Conjecture~4.8]{SYZ}.

\begin{corollary}[Sun--Yang--Zuo conjecture]\label{cor:syz}
Let $k$ be a perfect field of characteristic $p>2$, and
let $\lambda\in k\setminus\{0,1\}$. Then the following
diagram commutes.
\[
  \begin{tikzcd}
    E_\lambda \arrow[r,"{[p]}"] \arrow[d,"\pi"'] &
    E_\lambda \arrow[d,"\pi"] \\
    X \arrow[r,"\varphi_{\lambda,p}"'] & X
  \end{tikzcd}
\]
\end{corollary}

\begin{remark}
 The validity of the corollary does not directly imply that the theorem holds: the corollary is equivalent to assert equality of two fractions, while Theorem~\ref{thm_det} identifies the denominator determinants
up to an explicit nonzero factor.
The proof of \cite[Proposition~A.9]{SYZ} applies the determinant
identity at $(\lambda^{-1},a^{-1})$ to obtain the corresponding
relation between the numerator determinants and then substitutes
both relations into the formulas for the two self-maps.
\end{remark}

\subsection*{Historical remarks}

In \cite{LSW}, an approach to the conjecture is given by
lifting the Higgs--de Rham flow to a moduli space on the associated
elliptic curve. Lam and Litt credit
\cite{LSW} with an earlier proof and give their own proof using middle
convolution \cite[Theorem~1.11 and \S5]{LL}. 

Our proof follows the reduction of Sun, Yang and Zuo to a determinant identity. It establishes the determinant identity in \cite[Conjecture~A.1]{SYZ} by polynomial division and uses the reduction in \cite[Proposition~A.9]{SYZ}.

\subsection*{Acknowledgements.}
The author is supported by the Fundamental Research Funds for the Central Universities, the National Natural Science Foundation of China under Grant No. 12671021 and No. 12331002, and the Science Fund for Global Challenges and Sustainability under Grant No. W2441003.

\section{Proof of the determinant identity}

In this section, we retain the notation from Theorem~\ref{thm_det}.

\subsection{The key linear map $\Phi$ and its determinant}
\label{subsec_target-bases}

Set $R=k[t]$, $U=R[x]_{\le m}$, and $W=R[x]_{<m}$. 
Denote the $R$-linear map
\begin{equation}\label{eq_common-map}
  \begin{aligned}
    \Phi\colon U\oplus W\oplus W
      &\longrightarrow R[x]_{\le 3m},\\
    (\varphi,\psi,\tau)&\longmapsto g^m \varphi-x^{m+1}\psi-(x^p-t)\tau.
  \end{aligned}
\end{equation}
Both modules are free of rank $3m+1$. On the source we use
the ordered basis
\[
  \mathcal B_s=
  \Bigl(
    (x^i,0,0)_{i=0}^{m},\
    (0,x^i,0)_{i=0}^{m-1},\
    (0,0,x^i)_{i=0}^{m-1}
  \Bigr).
\]
On the target, division by $x^p-t$ and by $g^m$ gives,
respectively, the ordered bases
\begin{align*}
  \mathcal B_t
    &=\bigl(1,x,\ldots,x^{2m},\
            x^p-t,\ldots,x^{m-1}(x^p-t)\bigr),\\
  \mathcal B'_t
    &=\bigl(g^m,xg^m,\ldots,x^mg^m,\
            1,x,\ldots,x^{2m-1}\bigr).
\end{align*}
Let $M_1$ and $M_2$ be the matrices of $\Phi$ with
source basis $\mathcal B_s$ and target bases $\mathcal B_t$
and $\mathcal B'_t$, respectively. I.e., 
\[
  \Phi(\mathcal B_s)=\mathcal B_t M_1=\mathcal B'_t M_2,
\]
where $\Phi(\mathcal B_s)$ is the row of images of the
source basis vectors. Let $T$ be the change-of-coordinates
matrix defined by $\mathcal B'_t=\mathcal B_t T$. Then
\begin{equation}\label{eq_key}
  M_1=TM_2.
\end{equation}

\begin{lemma}\label{lem_target-bases}
The matrices $T$ and $M_1$ satisfy
\[
  \textup{(i)}\quad \det T=1,
  \qquad
  \textup{(ii)}\quad \det M_1=\det B.
\]
\end{lemma}

\begin{proof}
\emph{(i) Change of basis.}
Let $\mathcal B_t^{st}=(1,x,\ldots,x^{3m})$ be the standard
monomial basis. Reordering $\mathcal B'_t$ by increasing
degree gives
\[
  {\mathcal B'_t}^+
  =\bigl(1,x,\ldots,x^{2m-1},\
         g^m,xg^m,\ldots,x^mg^m\bigr).
\]
Define $T_1$, $T_2$, and $T_3$ by
\[
  \mathcal B_t=\mathcal B_t^{st}T_1,
  \qquad
  {\mathcal B'_t}^+=\mathcal B_t^{st}T_2,
  \qquad
  \mathcal B'_t={\mathcal B'_t}^+ T_3.
\]
The vectors of $\mathcal B_t$ and ${\mathcal B'_t}^+$
are monic polynomials of successive degrees $0,\ldots,3m$.
Their coordinate matrices $T_{1}$ and
$T_{2}$ are therefore upper triangular with
diagonal entries $1$. The permutation matrix $T_3$
interchanges blocks of lengths $2m$ and $m+1$, so
\[
  \det T_3=(-1)^{2m(m+1)}=1.
\]
Comparing $\mathcal B'_t = {\mathcal B'_t}^+ T_3 = \mathcal B_t^{st}T_2T_3$ 
with $\mathcal B'_t=\mathcal B_t T = \mathcal B_t^{st}T_1 T$ gives $T_2T_3 = T_1 T$. Hence $T=T_1^{-1}T_2T_3$ has determinant $1$.

\emph{(ii) Reduction modulo $x^p-t$.}
Divide $g^m\varphi$ by $x^p-t$. The quotient has degree at most
$3m-p=m-1$, and the remainder has degree at most
$p-1=2m$. Splitting the remainder into terms of degree
at most $m$ and terms of higher degree gives
\[
  g^m \varphi=\varphi_0+x^{m+1}\psi_0+(x^p-t)\tau_0,
  \qquad \varphi_0\in U,\quad \psi_0,\tau_0\in W.
\]
It follows that
\[
  \Phi(\varphi,\psi,\tau)
    =\varphi_0+x^{m+1}(\psi_0-\psi)+(x^p-t)(\tau_0-\tau).
\]
The three blocks of $\mathcal B_t$-coordinates are the
coefficient vectors of $\varphi_0$, $\psi_0-\psi$, and $\tau_0-\tau$.

To find the matrix of $\varphi\mapsto \varphi_0$, take $\varphi=x^s$,
where $0\le s\le m$. The definition of $\gamma_n$ gives
\[
  g^mx^s=\sum_{i=0}^{2m}\gamma_{m+i}x^{i+s}.
\]
Since $i+s\le 3m<2p$, only terms of degrees $r$ and
$p+r$ can contribute to the coefficient of $x^r$ in
the remainder. For $0\le r,s\le m$, there are two cases.

\begin{itemize}
  \item If $r\ge s$, then $0\le r-s\le m$ and
    $p+r-s>2m$. Only $i=r-s$ occurs in the sum, giving
    the coefficient $\gamma_{m+r-s}$.
  \item If $r<s$, then $r-s<0$ and
    $m+1\le p+r-s\le 2m$. Only $i=p+r-s$ occurs.
    Since $x^{p+r}\equiv tx^r\pmod{x^p-t}$, this
    gives the coefficient $t\gamma_{m+p+r-s}$.
\end{itemize}

These are the coefficients of $\varphi_0$, since $r\le m$,
and they agree with the entries $B_{rs}$. Thus the matrix
of $\varphi\mapsto \varphi_0$ in the ordered monomial basis of $U$
is $B$, and
\[
  M_1=
  \begin{pmatrix}
    B & 0 & 0\\
    * & -I_m & 0\\
    * & 0 & -I_m
  \end{pmatrix}.
\]
Taking determinants gives
\begin{equation}\label{eq_phi-b}
  \det M_1=\det B\,(-1)^{2m}=\det B.
\end{equation}
\end{proof}

\subsection{Reduction modulo \texorpdfstring{$g^m$}{g to the m}}
\label{subsec_remainder-map}

In the basis $\mathcal B'_t$, the last $2m$ coordinates
record the remainder modulo $g^m$. For a monic polynomial
$h$, write $\operatorname{rem}_h(f)$ for the remainder
of $f$ upon division by $h$. Define
\[
  \begin{aligned}
    \Psi\colon W\oplus W&\longrightarrow R[x]_{<2m},\\
    (\psi,\tau)&\longmapsto
      \operatorname{rem}_{g^m}\bigl(x^{m+1}\psi+(x^p-t)\tau\bigr).
  \end{aligned}
\]
For a map $F$ among $\Psi$ and the maps below, write
\[[F]\]
for its matrix in the monomial bases ordered by increasing degree, with
direct-sum factors taken in the displayed order.

\begin{lemma}\label{lem_det-S}
The lower-right $2m\times 2m$ block of $M_2$ is
$-[\Psi]$, and
\[
  \det M_2=\det[\Psi].
\]
\end{lemma}

\begin{proof}
For $0\le i\le m$, the image
$\Phi(x^i,0,0)=x^ig^m$ is the corresponding basis vector
in the first block of $\mathcal B'_t$. The first $m+1$
columns of $M_2$ therefore give the blocks $I_{m+1}$
and $0$.

For the remaining columns of $M_2$, we have
\[
  \operatorname{rem}_{g^m}\bigl(\Phi(0,\psi,\tau)\bigr)
    =-\operatorname{rem}_{g^m}
      \bigl(x^{m+1}\psi+(x^p-t)\tau\bigr)
    =-\Psi(\psi,\tau).
\]
The source vectors for these columns are ordered as in
the source basis of $[\Psi]$. Thus the lower-right block
of $M_2$ is $-[\Psi]$, and
\[
  M_2=
  \begin{pmatrix}
    I_{m+1} & *\\
    0 & -[\Psi]
  \end{pmatrix}.
\]
Since $[\Psi]$ has size $2m\times 2m$, we obtain
\begin{equation}\label{eq_b-s}
  \det M_2=\det(-[\Psi])=(-1)^{2m}\det[\Psi]=\det[\Psi].
\end{equation}
\end{proof}

\subsection{Local remainders}
\label{subsec_local-remainders}

The factorization $g^m=(x-1)^m(x-\lambda)^m$ lets us read
the output of $\Psi$ through its remainders at $1$ and
$\lambda$. For $z\in\{1,\lambda\}$, set 
\[ \rho_z \colon  R[x]  \to W, \qquad 
  f\mapsto \operatorname{rem}_{(x-z)^m}(f),\]
\[\Gamma\colon R[x]_{<2m} \to W\oplus W, \qquad   f\mapsto \bigl(\rho_1(f),\rho_\lambda(f)\bigr).
\]
and
\[\Lambda_z \colon W\to W, \qquad  \Lambda_z(f)=\rho_z(x^mf).\]
We obtain the composite
\[
  W\oplus W
    \xrightarrow{\,\Psi \,} R[x]_{<2m}
    \xrightarrow{\,\Gamma\,} W\oplus W
    \xrightarrow{\,\Lambda_1\oplus \Lambda_\lambda\,} W\oplus W,
\]
whose matrix is
\begin{equation}\label{eq_K-product}
  M'_2=\begin{pmatrix} [\Lambda_1] & \\ &[\Lambda_\lambda] \end{pmatrix} \cdot  [\Gamma]\cdot [\Psi].
\end{equation}
We also set 
\begin{equation} \label{eq_Z-definition}
M''_2=(\lambda^p-t)[\Lambda_\lambda]-\lambda^p(1-t)[\Lambda_1]
\end{equation}

\begin{lemma}\label{lem_local-remainders}
We have 
\begin{enumerate}
    \item $\det[\Gamma]=(\lambda-1)^{m^2}$,
  $\det[\Lambda_z]=z^{m^2}$ for $z\in\{1,\lambda\})$, and 
  \item $\det M''_2=\det M'_2
    =\lambda^{m^2}(\lambda-1)^{m^2}\det[\Psi]$.
\end{enumerate}
\end{lemma}

\begin{proof}
The polynomials $(x-1)^m$ and $(x-\lambda)^m$ are
relatively prime in $k[x]$, so they generate the unit ideal
in $R[x]$. The Chinese remainder theorem gives the
isomorphism
\[
  R[x]/(g^m)
    \longrightarrow
    R[x]/\bigl((x-1)^m\bigr)
       \oplus R[x]/\bigl((x-\lambda)^m\bigr).
\]
Representing each class by its unique remainder identifies
this isomorphism with $\Gamma$.

To compute $\det[\Gamma]$, temporarily use the ordered
source basis
\[
  \bigl((x-1)^j\bigr)_{j=0}^{m-1},\qquad
  \bigl((x-1)^m(x-\lambda)^j\bigr)_{j=0}^{m-1},
\]
and the ordered basis $1,x-z,\ldots,(x-z)^{m-1}$ in
the target component indexed by $z$. Each change of
basis from the monomial basis has determinant $1$,
since the new polynomials are monic and have successive
degrees in the corresponding space.

In these bases the first component of $\Gamma$ gives
the block row $(I_m,0)$. For the lower-right block, put
$y=x-\lambda$. The remainder of
$(x-1)^m(x-\lambda)^j$ modulo $(x-\lambda)^m$ is
\[
  \sum_{a=0}^{m-1-j}
    \binom{m}{a}(\lambda-1)^{m-a}y^{a+j},
\]
obtained by expanding $(\lambda-1+y)^m y^j$ and
discarding terms of degree at least $m$ in $y$.
The coefficient of $y^j$ is $(\lambda-1)^m$,
and all coefficients of lower degree vanish. Thus the
lower-right block is lower triangular with diagonal
entries $(\lambda-1)^m$. Taking determinants gives
\[
  \det[\Gamma]=\bigl((\lambda-1)^m\bigr)^m
    =(\lambda-1)^{m^2}.
\]

For $\Lambda_z$, use the basis $1,y,\ldots,y^{m-1}$ of $W$,
where $y=x-z$. Expanding $(z+y)^m y^j$ and discarding
terms of degree at least $m$ in $y$ gives
\[
  \Lambda_z(y^j)
    =\sum_{a=0}^{m-1-j}\binom{m}{a}z^{m-a}y^{a+j}.
\]
The resulting matrix is lower triangular with diagonal
entries $z^m$. Using the same temporary basis on the
source and target preserves the determinant, so
\[
  \det[\Lambda_z]=(z^m)^m=z^{m^2}.
\]

In the monomial bases, the component indexed by $z$
of $(\Lambda_1\oplus \Lambda_\lambda)\circ\Gamma\circ \Psi$ sends
$(\psi,\tau)$ to
\[
  \Lambda_z\bigl(\rho_z(\Psi(\psi,\tau))\bigr)
    =\rho_z\bigl(x^m\Psi(\psi,\tau)\bigr).
\]
By definition, $\Psi(\psi,\tau)$ is congruent to
$x^{m+1}\psi+(x^p-t)\tau$ modulo $g^m$, hence also modulo
$(x-z)^m$. Since $x^p-z^p=(x-z)^p$ in characteristic
$p$, we have $x^p\equiv z^p\pmod{(x-z)^m}$. Using
$p=2m+1$ and $\psi\in W$, we obtain
\begin{align*}
  \rho_z\bigl(x^m\Psi(\psi,\tau)\bigr)
    &=\rho_z\bigl(x^p\psi+x^m(x^p-t)\tau\bigr)\\
    &=z^p\psi+(z^p-t)\rho_z(x^m\tau)\\
    &=z^p\psi+(z^p-t)\Lambda_z(\tau).
\end{align*}
Taking $z=1$ and $z=\lambda$ in that order gives
\[
  M'_2=
  \begin{pmatrix}
    I_m & (1-t)[\Lambda_1]\\
    \lambda^p I_m & (\lambda^p-t)[\Lambda_\lambda]
  \end{pmatrix}.
\]
Subtracting $\lambda^p$ times the first block row from
the second preserves the determinant and gives
\[
  \begin{pmatrix}
    I_m & (1-t)[\Lambda_1]\\
    0 & M''_2
  \end{pmatrix},
\]
where $M''_2$ is given by \eqref{eq_Z-definition}. Hence
$\det M'_2=\det M''_2$. The product formula
\eqref{eq_K-product} gives 
\begin{align*}
  \det M'_2
    &=\det[\Lambda_1]\,\det[\Lambda_\lambda]\,\det[\Gamma]\,\det[\Psi]\\
    &=\lambda^{m^2}(\lambda-1)^{m^2}\det[\Psi]. \qedhere
\end{align*}
\end{proof}

\subsection{Comparison with \texorpdfstring{$A$}{A}}
\label{subsec_comparison-A}

To compare $M''_2$ with $A$, we separate the row and column
factors in the entries of $[\Lambda_z]$. Set
\[
  D_L=\operatorname{diag}(d_0,\ldots,d_{m-1}),
  \qquad
  D_R=\operatorname{diag}(e_0,\ldots,e_{m-1}),
\]
where
\[
  d_r=\frac{(-1)^{m-1-r}}{r!(m-1-r)!},
  \qquad
  e_j=\frac{(m+j)!}{j!}
  \qquad (0\le r,j<m).
\]
All factorial arguments are less than $p$, so $D_L$
and $D_R$ are invertible.

\begin{lemma}\label{lem_comparison-A}
With $C$ and $c_p$ as in Theorem~\ref{thm_det}, we have
\[
  M''_2=-D_LAD_R,\qquad [\Lambda_1]=D_LCD_R.
\]
Consequently, $c_p\in\mathbb F_p^\times$ and
$\det A=c_p\det M''_2$.
\end{lemma}

\begin{proof}
For $0\le r,j<m$, the $(r,j)$-entry of $[\Lambda_z]$
is the coefficient of $x^r$ in $\rho_z(x^{m+j})$.
Expanding $x^{m+j}=(z+(x-z))^{m+j}$ and discarding terms
divisible by $(x-z)^m$ gives
\[
  \rho_z(x^{m+j})
    =\sum_{\ell=0}^{m-1}
       \binom{m+j}{\ell}z^{m+j-\ell}(x-z)^\ell.
\]
Put $n=m+j-r$ and $b=m-1-r$. Extracting the
coefficient of $x^r$ gives
\begin{align*}
  [\Lambda_z]_{rj}
    &=\sum_{\ell=r}^{m-1}
      \binom{m+j}{\ell}\binom{\ell}{r}
      z^{m+j-\ell}(-z)^{\ell-r}\\
    &=z^n\binom{m+j}{r}
      \sum_{q=0}^{b}(-1)^q\binom{n}{q}.
\end{align*}
In the second line we set $q=\ell-r$ and used
\[
  \binom{m+j}{\ell}\binom{\ell}{r}
    =\binom{m+j}{r}\binom{m+j-r}{\ell-r}.
\]
Pascal's identity gives
\begin{align*}
  \sum_{q=0}^{b}(-1)^q\binom{n}{q}
    &=\sum_{q=0}^{b}(-1)^q\binom{n-1}{q}
      -\sum_{q=0}^{b-1}(-1)^q\binom{n-1}{q}\\
    &=(-1)^b\binom{n-1}{b},
\end{align*}
where the second sum is empty if $b=0$.
Since $n-1-b=j$, we have
\[
  \binom{m+j}{r}\binom{n-1}{b}
    =\frac{(m+j)!}{r!\,n!}\frac{(n-1)!}{b!\,j!}
    =\frac{1}{r!\,b!}\frac{1}{n}\frac{(m+j)!}{j!}.
\]
All factorial arguments are less than $p$, and
$1\le n\le p-2$. Substituting $d_r$ and $e_j$ gives
\[
  [\Lambda_z]_{rj}=d_r\,\frac{z^n}{n}\,e_j,
  \qquad n=m+j-r.
\]

By \eqref{eq_Z-definition}, this entry formula gives
\begin{align*}
  (M''_2)_{rj}
    &=(\lambda^p-t)d_r\frac{\lambda^n}{n}e_j
      -\lambda^p(1-t)d_r\frac{1}{n}e_j\\
    &=-d_r\,\frac{\lambda^p(1-t)-(\lambda^p-t)\lambda^n}{n}\,e_j\\
    &=-d_r A_{rj}e_j.
\end{align*}
Thus $M''_2=-D_LAD_R$. Taking $z=1$ in the same entry
formula gives $[\Lambda_1]_{rj}=d_r C_{rj}e_j$, hence
$[\Lambda_1]=D_LCD_R$.

By Lemma~\ref{lem_local-remainders}, $\det[\Lambda_1]=1$.
Taking determinants in the second factorization gives
\[
  1=\det D_L\,\det C\,\det D_R,
  \qquad
  \det C=(\det D_L\,\det D_R)^{-1}.
\]
In particular, $\det C\ne0$. Since the entries of $C$
lie in $\mathbb F_p$, we have
$c_p=(-1)^m\det C\in\mathbb F_p^\times$.
Taking determinants in the first factorization gives
\[
  \det M''_2=(-1)^m\det D_L\,\det A\,\det D_R.
\]
Multiplying by $c_p=(-1)^m\det C$ and substituting
$\det C=(\det D_L\,\det D_R)^{-1}$ gives
$\det A=c_p\det M''_2$.
\end{proof}

\begin{proof}[Proof of Theorem~\ref{thm_det}]
By \eqref{eq_key} and Lemma~\ref{lem_target-bases},
$\det M_2= \det M_1=\det B$, since $\det T=1$.
Lemmas~\ref{lem_det-S}, \ref{lem_local-remainders},
and~\ref{lem_comparison-A} now give
\begin{align*}
  \det A
    &=c_p\det M''_2\\
    &=c_p\,\lambda^{m^2}(\lambda-1)^{m^2}\det[\Psi]\\
    &=c_p\,\lambda^{m^2}(\lambda-1)^{m^2}\det M_2\\
    &=c_p\,\lambda^{m^2}(\lambda-1)^{m^2}\det B.
\end{align*}
Lemma~\ref{lem_comparison-A} also gives
$c_p\in\mathbb F_p^\times$.
\end{proof}

\end{document}